\documentclass[]{amsart}
\usepackage{times,mathrsfs,color,comment}
\usepackage{amssymb,amsmath,bbm}

\numberwithin{equation}{section}

\newtheorem{thm}{Theorem}[section]
\newtheorem{cor}[thm]{Corollary}
\newtheorem{lem}[thm]{Lemma}
\newtheorem{ex}[thm]{Example}
\newtheorem{pro}[thm]{Proposition}

\newtheorem{rmk}[thm]{Remark}

\newcommand{\re}{{\rm Re\,}}

\newcommand{\D}{{\mathbb D}}

\begin{document}

\title[]
{Kernel-induced distance and its applications to Composition operators on Large Bergman spaces}

\author[I. Park]{Inyoung Park}
\address{(Park) Institute of Mathematical Sciences, Ewha Womans University, Seoul 03760, KOREA}
\email{iypark26@gmail.com}

\thanks{The author was supported by NRF (No. 2019R1A6A1A11051177) of Korea.}


\begin{abstract}
In this paper, we obtain a complete characterization for the compact difference of two composition operators acting on Bergman spaces with a rapidly decreasing weight $\omega=e^{-\eta}$, $\Delta\eta>0$. In addition, we provide simple inducing maps which support our main result. We also study the topological path connected component of the space of all bounded composition operators on $A^2(\omega)$ endowed with the Hilbert-Schmidt norm topology.
\end{abstract}

\maketitle

\section{Introduction}
Let $S(\D)$ be the space of holomorphic self-maps of the unit disk $\D$. Given $\varphi$ in $S(\D)$, the composition operator $C_\varphi$ is defined by $C_\varphi f:=f\circ\varphi$ for all $f$ belonging to the holomorphic function spaces $H(\D)$. For the integrable radial function $\omega$, let $L^p(\omega dA)$ be the space of all measurable functions $f$ on $\D$ such that
\begin{align*}
\|f\|^p_{p}:=\int_\D |f(z)\omega(z)|^pdA(z)<\infty,\quad0<p<\infty,
\end{align*}
where $dA(z)$ is the normalized area measure on $\D$. We put $\|\cdot\|:=\|\cdot\|_2$ for simplicity. Especially, we denote $A^p(\omega):=L^p(\omega dA)\cap H(\D)$ and we use the notation $A^p_\alpha(\D)$ when $\omega(z)=(1-|z|)^\alpha$, $\alpha>-1$. Throughout this paper, we consider positive radial weights of the form $\omega=e^{-\eta}$ where $\eta$ is a strictly increasing radial function and $\Delta\eta>0$ on $\D$. Now, it is said that the weight $\omega$ belongs to the class $\mathcal{W}$ if we can choose a differentiable radial function $\tau$ as
\begin{align}\label{taudef}
(\Delta\eta(z))^{-\frac{1}{2}}\asymp\tau(r)
\end{align}
which $\lim_{r\rightarrow1^-}\tau(r)=0$ and $\lim_{r\rightarrow1^-}\tau'(r)=0$. We remark that $\mathcal{W}$ is a subclass of the so-called fast weights considered by Kriete and MacCluer in \cite{CM,KM}. Especially, $\eta(z)=\frac{1}{1-|z|}$ is contained in $\mathcal{W}$ as a typical example. Readers can refer to \cite{HLS,PP} to study other admissible weights belonging to $\mathcal{W}$. On the other hand, we note that the standard weights $\eta(r)=-\alpha\log(1-r)$, $\alpha>-1$, are not contained in the class $\mathcal{W}$ since we can not find $\tau$ such that $\tau'(r)$ converge to $0$ when $r\rightarrow1^-$. That means we do not consider the standard weighted Bergman spaces in this paper.\\
\indent It is well known that all composition operators are bounded in the standard weighted Bergman spaces $A^p_\alpha(\D)$ by the Littlewood subordination principle, but it is not true in Bergman spaces with fast weights any more. Much studies on the theory of composition operators over the large Bergman spaces has been established in \cite{CM,KM}. When $\omega$ belongs to the class $\mathcal{W}$, it has known that $C_\varphi$ is bounded on $A^2(\omega)$  if and only if
\begin{align}\label{bddwithoutU}
\sup_{z\in\D}\frac{\omega(z)}{\omega(\varphi(z))}<\infty.
\end{align}
You can also refer to the proof in \cite[Theorem 3.2]{P} for the nessecity. As a consequence of the Carleson measure theorem, we observe that the condition (\ref{bddwithoutU}) still holds for the boundedness of $C_\varphi$ on $A^p(\omega)$ for all range of $p$ (See Remark \ref{bddp}). For the compactness, Kriete and MacCluer gave the estimate of the essential norm of $C_\varphi$ on $A^2(\omega)$:
\begin{align*}
\|C_\varphi\|_e\approx\limsup_{|z|\rightarrow1}\frac{\omega(z)}{\omega(\varphi(z))}.
\end{align*}

In the case of classical holomorphic function spaces such as the Hardy and the standard weighted Bergman space over the disk or the ball, many characterizations for the compactness of $C_\varphi-C_\psi$ have been developed over the past decades involving the pseudo-hyperbolic distance, but no result has been given in the large Bergman space; see for example \cite{CCKY, KW, M, MOZ}. In order to state our main result, we introduce the following distance induced by function $\tau$ associated with the weight $\omega\in\mathcal{W}$ as in (\ref{taudef}):
\begin{align*}
\rho_{\tau,\varphi,\psi}(z):=\rho_\tau(\varphi(z),\psi(z))=1-e^{-d_\tau(\varphi(z),\psi(z))}
\end{align*}
where
\begin{align}\label{dist}
d_\tau(z,w):=\inf_\gamma\int^1_0\frac{|\gamma'(t)|}{\tau(\gamma(t))}dt,
\end{align}
where the infimum is taken over all piecewise smooth curves $\gamma$ connecting $z$ and $w$. One may check that $\rho_{\tau}(z,w)$ can be a distance in \cite[Lemma 3.3]{P1}. In Section 4, we give the necessary conditions for the boundedness and compactness of difference of composition operators in terms of the distance involving $\tau$ function:
\begin{thm}\label{k-thm}
Let $\omega\in\mathcal{W}$. If $C_\varphi-C_\psi$ is bounded (compact, resp.) on $A^p(\omega)$ for $0<p<\infty$ then
\begin{align}\label{functheor}
\sup_{z\in\D}\rho_{\tau,\varphi,\psi}(z)^2\left(\frac{\omega(z)}{\omega(\varphi(z))}+\frac{\omega(z)}{\omega(\psi(z))}\right)<\infty \ \  (\longrightarrow0\ \ as \ \ |z|\rightarrow1^-, resp.)
\end{align}
\end{thm}
For the compact difference of composition operators on $A^p(\omega)$, we obtain the following equivalent condition:
\begin{thm}\label{k-thm1}
Let $\omega\in\mathcal{W}$ and $C_\varphi$, $C_\psi$ be bounded on $A^p(\omega)$ for $0<p<\infty$. Then the following statements are equivalent: for $0<p<\infty$,
\begin{itemize}
  \item [(1)]$C_\varphi-C_\psi$ is compact on $A^p(\omega)$.
  \item [(2)]$\lim_{|z|\rightarrow1^-}\rho_{\tau,\varphi,\psi}(z)^2\left(\frac{\omega(z)}{\omega(\varphi(z))}+\frac{\omega(z)}{\omega(\psi(z))}\right)=0$.
  \item [(3)]$\rho_{\tau,\varphi,\psi} C_\varphi$ and $\rho_{\tau,\varphi,\psi} C_\psi$ are compact from $A^p(\omega)$ into $L^p(\omega dA)$.
\end{itemize}
\end{thm}
In fact, the current author couldn't prove that \eqref{functheor} is a lower bound for $\|C_\varphi-C_\psi\|_p$ in \cite{P1}, which makes it impossible to obtain the equivalent condition as presented in Theorem \ref{k-thm1}. But in this paper, we overcome the difficulty using the Skwarczy\`nski distance, which is defined by the reproducing kernel $\overline{K_z(w)}=K(z,w)$ in $A^2(\omega)$:
\begin{align*}
\mathcal{S}(z,w):=\sqrt{1-\frac{|K(z,w)|}{\|K_z\|\|K_w\|}},\quad z,w\in\D.
\end{align*}
It can be easily seen that $\mathcal{S}(z,w)$ is comparable to the pseudo-hyperbolic distance $\rho(z,w)=\big|\frac{z-w}{1-\overline{z}w}\big|$ in the case of classical weighted Bergman spaces and Hardy space. However, non-existence of the explicit form of the reproducing kernel in $A^2(\omega)$ makes it difficult to verify such a relation with $\mathcal{S}(z,w)$. In this paper, we give some inequality between $\rho_\tau(z,w)$ and $\mathcal{S}(z,w)$ in Theorem \ref{key1} and use it to obtain the norm estimate for the difference of the reproducing kernels in Theorem \ref{diffkernelnorm}. More information on the Skwarczy\'nski distance will be presented in Section 3 and other essential materials to be used for this paper will be stated in Section 2.\\ \indent In Section 5, we show that $C_\varphi-C_\psi$ is compact on $A^2(\omega)$ when the inducing maps $\varphi, \psi$ have good boundary behavior in the sense of higher-order data and order of contact. Using this result, we provide explicit analytic maps $\varphi, \psi$ which induce non-compact composition operators $C_\varphi$ and $C_\psi$ but $C_\varphi-C_\psi$ is compact on $A^2(\omega)$. In Section 6, we characterize a topological path component of $\mathcal{C}_{HS}(A^2(\omega))$, the space of all composition operators on $A^2(\omega)$ endowed with the topology induced by the Hilbert-Schmidt norm. The proofs of the following theorem are provided in Theorem \ref{mainresult5} and Theorem \ref{mainresult1}.
\begin{thm}\label{k-thm3}
Let $\omega\in\mathcal{W}$. Then the following statements are equivalent:
\begin{itemize}
  \item [(1)]$\rho_{\tau,\varphi,\psi}(z)\|K_{\varphi(z)}\|, \rho_{\tau,\varphi,\psi}(z)\|K_{\psi(z)}\|$ are belong to $L^2(\omega dA)$.
  \item [(2)]$C_\varphi-C_\psi$ is a Hilbert-Schmidt operator on $A^2(\omega)$.
  \item [(3)]$C_\varphi$ and $C_\psi$ lie in the same path component of $\mathcal{C}_{HS}(A^2(\omega))$.
\end{itemize}
\end{thm}

\bigskip

\textit{Constants.}
 In the rest of this paper, we use the notation
$X\lesssim Y$  or $Y\gtrsim X$ for nonnegative quantities $X$ and $Y$ to mean
$X\le CY$ for some inessential constant $C>0$. Similarly, we use the notation $X\approx Y$ if both $X\lesssim Y$ and $Y \lesssim X$ hold.

\bigskip
\section{Preliminary}
Assuming that $\tau'(r)\rightarrow0$ as $r\rightarrow1^-$, there exist constants $c_1, c_2>0$ independent of $z,w$ such that $\tau(z)< c_1(1-|z|)$ and
\begin{align}\label{comparable}
|\tau(z)-\tau(w)|\leq c_2|z-w|,\qquad z,w\in\D.
\end{align}
Throughout this paper, we denote
\begin{align*}
m_\tau:=\frac{\min(1,c_1^{-1},c_2^{-1})}{4}.
\end{align*}
\subsection{Radius functions and associated distance}
We let $D(z,r)$ be a Euclidean disk centred at $z$ with radius $r>0$. Using (\ref{comparable}), we obtain that for $0<\delta\leq m_\tau$,
\begin{align}\label{equiquan}
\frac{1}{2}\tau(z)\leq\tau(w)\leq2\tau(z)\quad\text{if}\quad w\in D(z,\delta\tau(z)),
\end{align}
and we use the notation $D(\delta\tau(z)):=D(z,\delta\tau(z))$ for simplicity. We set all $\delta$ appearing in the rest of our paper to meet the  conditions above. In \cite{AP,HLS}, when $\omega\in\mathcal{W}$, the authors gave the following useful inequality to estimate the reproducing kernel function of $A^2(\omega)$: for each positive number $M>0$, there exists a constant $C(M)>1$ such that
\begin{align}\label{estimate}
e^{-d_\tau(z,w)}\leq C\left(\frac{\min(\tau(z),\tau(w))}{|z-w|}\right)^M,\quad z\neq w\in\D.
\end{align}
The following lemma, together with (\ref{dist}), shows that $d_\tau(z,w)$ and $\frac{|z-w|}{\tau(z)}$ are comparable when $d_\tau(z,w)<R$ for some $R>0$.
\begin{lem}\label{setinclus}
Let $\omega\in\mathcal{W}$. If $d_\tau(z,w)<R$ for some $R>0$ then there is a constant $C_1>0$ depending only on constants $R$ and $M,C$ in (\ref{estimate}) such that
\begin{align*}
d_\tau(z,w)\geq C_1\frac{|z-w|}{\min(\tau(z),\tau(w))}.
\end{align*}
\end{lem}
\begin{proof}
For a given $z\in\D$ and any point $w\neq z$ with $d_\tau(z,w)<R$, from (\ref{estimate})
\begin{align*}
\frac{|z-w|}{\min(\tau(z),\tau(w))}<(Ce^R)^{\frac{1}{M}}.
\end{align*}
Thus, assuming that $|z|\geq|w|$, there exists a constant $0<s\leq1/\delta(Ce^R)^{\frac{1}{M}}$ satisfying
\begin{align}\label{distzw}
|z-w|=s\delta\tau(z).
\end{align}
We let $\gamma$ denote any curve satisfying that $\gamma(0)=z$, $\gamma(1)=w$ and $d_\tau(z,w)=\int^{1}_0\frac{|\gamma'(t)|}{\tau(\gamma(t))}dt$. Choose the minimum value $0<t_0\leq1$ such that
\begin{align}\label{distzw1}
|z-\gamma(t_0)|=\delta\tau(z).
\end{align}
For the case of $1<s\leq1/\delta(Ce^R)^{\frac{1}{M}}$ taken in (\ref{distzw}), by (\ref{equiquan}) and (\ref{distzw1}), we have
\begin{align}\label{s<1}
d_\tau(z,w)\geq\int^{t_0}_0\frac{|\gamma'(t)|}{\tau(\gamma(t))}dt\geq\frac{1}{2\tau(z)}\int^{t_0}_0|\gamma'(t)|dt&\geq\frac{\delta}{2}=\frac{|z-w|}{2s\tau(z)}\\
&\geq\frac{\delta}{2}(Ce^R)^{-\frac{1}{M}}\frac{|z-w|}{\tau(z)}.\nonumber
\end{align}
For $0<s\leq1$, we obtain $d_\tau(z,w)\geq\frac{|z-w|}{2\tau(z)}$ promptly from (\ref{s<1}). Thus, if we choose $C_1=\delta/2(Ce^R)^{-\frac{1}{M}}$ then we obtain the desired inequality.
\end{proof}

\subsection{Sub-mean value type inequalities}
The following inequalities play a crucial role in our proofs. In fact, their proofs are similar to the case of a doubling measure $\Delta\eta$ whose proof we can find in \cite[Lemma 19]{MMO}. For the proofs in our setting, you can refer to \cite{AP,HLS,O,PP}.
\begin{lem}\label{derivsubmean}
Let $\omega=e^{-\eta}$, where $\eta$ is a subharmonic function and $0<p<\infty$. Suppose the function $\tau$ satisfies properties (\ref{equiquan}) and  $\tau(z)^2\Delta\eta(z)\lesssim1$. Given $\delta>0$ satisfying (\ref{equiquan}) and $f\in H(\D)$,
\begin{itemize}
  \item[(1)]$|f(z)e^{-\eta(z)}|^p\lesssim\frac{1}{\tau(z)^2}\int_{D(\delta\tau(z))}|f(\xi)e^{-\eta(\xi)}|^pdA(\xi)$,
  \item[(2)]$|f'(z)e^{-\eta(z)}|^p\lesssim\frac{1}{\tau(z)^{2+p}}\int_{D(\delta\tau(z))}|f(\xi) e^{-\eta(\xi)}|^p\,dA(\xi)$.
\end{itemize}
\end{lem}
Now, we prove the following inequality which gives the upper estimate for the difference in two variables $z, w$ of a function with $|z-w|<\delta\min(\tau(z),\tau(w))$.
\begin{lem}\label{mvtforinterg}
Let $\omega=e^{-\eta}\in\mathcal{W}$ and $0<p<\infty$. Then for $f\in H(\D)$ and $w\in D(\delta\tau(z))$ with $|z|\geq|w|$,
\begin{align*}
|f(z)-f(w)|^pe^{-p\eta(z)}\lesssim\frac{\rho_\tau(z,w)^p}{\tau(z)^{2}}\int_{D(6\delta\tau(z))}|f(\xi) e^{-\eta(\xi)}|^p\,dA(\xi).
\end{align*}
The constant suppressed depends on $C,M,\delta$ and $p$.
\end{lem}
\begin{proof}
By the fundamental theorem of integration, for $w\in D(\delta\tau(z))$,
\begin{align}\label{zt}
|f(z)-f(w)|&=\left|\int_{0}^{1}f'(zt+w(1-t))(z-w)dt\right|\nonumber\\&\leq|z-w|\sup_{t\in[0,1]}|f'(zt+w(1-t))|.
\end{align}
Denote $z_t:=zt+w(1-t)$, $0\leq t\leq1$. By (\ref{equiquan}), for $\xi\in D(\delta\tau(z_t))$, we have
\begin{align}\label{setinclus1}
|\xi-z_i|&\leq|\xi-z_t|+|z_t-z_i|\nonumber\\&\leq\delta\tau(z_t)+\delta\tau(z)\leq3\delta\tau(z)\leq6\delta\tau(z_i)\quad\text{for}\quad i=0,1.
\end{align}
Since $|z_t|\leq\max(|z|,|w|)=|z|$ by our assumption, (2) of Lemma \ref{derivsubmean} and (\ref{setinclus1}) yield
\begin{align*}
|f'(z_t)|^p&\lesssim\frac{e^{p\eta(z_t)}}{\tau(z_t)^{2+p}}\int_{D(\delta\tau(z_t))}|f(\xi) e^{-\eta(\xi)}|^p\,dA(\xi)\\
&\lesssim\frac{e^{p\eta(z)}}{\tau(z)^{2+p}}\int_{D(6\delta\tau(z))}|f(\xi) e^{-\eta(\xi)}|^p\,dA(\xi).
\end{align*}
Applying the inequality above to (\ref{zt}), there is a constant $C:=C(\delta,p,M)>0$ by Lemma \ref{setinclus} such that
\begin{align*}
|f(z)-f(w)|^pe^{-p\eta(z)}\leq C\frac{d_\tau(z,w)^p}{\tau(z)^{2}}\int_{D(6\delta\tau(z))}|f(\xi) e^{-\eta(\xi)}|^p\,dA(\xi).
\end{align*}
Since $e^{-2\delta}x\leq1-e^{-x}$ for $0\leq x<2\delta<1$, we have
\begin{align*}
d_\tau(z,w)\leq e^{2\delta}\rho_{\tau}(z,w).
\end{align*}
Hence, we complete the proof.
\end{proof}

\subsection{Carleson measure theorem}
A positive Borel measure $\mu$ in $\D$ is called a (vanishing) Carleson measure for $A^p(\omega)$ if the embedding $A^p(\omega)\subset L^p(\omega d\mu)$ is (compact) continuous where
\begin{align*}
L^p(\omega d\mu):=\left\{f\in\mathcal{M}(\D)\big|\int_\D|f(z)\omega(z)|^pd\mu(z)<\infty\right\}
\end{align*}
and $\mathcal{M}(\D)$ is a set of $\mu$-measurable functions on $\D$. Now, we introduce Carleson measure theorem in our setting, which can be found in \cite{LR, O, PP} for example.
\begin{thm}[Carleson measure theorem]\label{embedding}
Let $\omega\in\mathcal{W}$ and $\mu$ be a positive Borel measure on $\D$. Then, for $0<p<\infty$, we have
\begin{itemize}
\item[(1)]The embedding $I:A^p(\omega)\rightarrow L^p(\omega d\mu)$ is bounded if and only if for a small $\delta\in(0,m_\tau)$, we have $\sup_{z\in\D}\frac{\mu(D(\delta\tau(z)))}{\tau(z)^{2}}<\infty$.
\item[(2)]The embedding $I:A^p(\omega)\rightarrow L^p(\omega d\mu)$ is compact if and only if for a small $\delta\in(0,m_\tau)$, we have $\lim_{|z|\rightarrow1}\frac{\mu(D(\delta\tau(z)))}{\tau(z)^{2}}=0$.
\end{itemize}
\end{thm}
By the measure theoretic change of variables, we have
\begin{align*}
\|uC_\varphi f\|_{p}^p=\int_\D|u(f\circ\varphi)\omega|^p dA=\int_\D|f\omega|^p d\mu_{u,\varphi,p}
\end{align*}
where
\begin{align*}
\mu_{u,\varphi,p}(E):=\omega^{-p}[|u|^p\omega^{p}dA]\circ\varphi^{-1}(E)=\int_{\varphi^{-1}(E)}|u(z)|^p\frac{\omega(z)^{p}}{\omega(\varphi(z))^{p}}dA(z)
\end{align*}
for any measurable subsets $E$ of $\D$. Denote the norm of operator $T$ acting on $A^p(\omega)$ by $\|T\|_{A^p(\omega)}$. It is well known that
\begin{align}\label{opnorm}
\|uC_\varphi\|^p_{A^p(\omega)}\approx\sup_{z\in\D}\frac{\mu_{u,\varphi,p}(D(\delta\tau(z)))}{\tau(z)^{2}}.
\end{align}
Since the weight $\omega^{2/p}$ still belongs to $\mathcal{W}$ and the same $\tau$ function of $A^p(\omega^{2/p})$ can be chosen with that of $A^2(\omega)$, we have
\begin{align}\label{bddindp}
\|C_\varphi\|^p_{A^p(\omega)}\approx\sup_{z\in\D}\frac{\omega^{-p}[\omega^pdA]\circ\varphi^{-1}(D(\delta\tau(z)))}{\tau(z)^{2}}\approx\|C_\varphi\|^2_{A^2(\omega^{p/2})}
\end{align}
The constant suppressed in the second estimate above depends on $p$ and $\delta$. As a consequence of (\ref{bddindp}), we conclude the following Remark:
\begin{rmk}\label{bddp}
Let $\omega\in\mathcal{W}$, $0<p<\infty$. Then $C_\varphi$ is bounded on $A^p(\omega)$ if and only if
\begin{align*}
\sup_{z\in\D}\frac{\omega(z)}{\omega(\varphi(z))}<\infty.
\end{align*}
Therefore, we conclude that $C_\varphi$ is bounded on $A^p(\omega)$ for all $0<p<\infty$ if it is bounded on $A^p(\omega)$ for some $0<p<\infty$.
\end{rmk}
\begin{proof}
From (\ref{bddindp}), the boundedness of $C_\varphi$ acting on $A^p(\omega)$ is equivalent to the boundedness of $C_\varphi$ acting on $A^2(\omega^{p/2})$. Thus, it is clear that the first assertion holds by (\ref{bddwithoutU}).
\end{proof}

\section{Kernel-induced pseudodistance}
Given $\omega\in\mathcal{W}$, the space $A^2(\omega)$ is a closed subspace of $L^2(\omega dA)$ with inner product
\begin{align*}
\langle f,g\rangle_\omega=\int_\D f\overline{g}\omega^2 dA,\quad f,g\in A^2(\omega).
\end{align*}
As is well known, the reproducing kernel of the Bergman space $A^2(\omega)$ is defined by
\begin{align}\label{definition}
K(z,w)=\overline{K_z(w)}=\sum_{k=0}^\infty e_k(z)\overline{e_k(w)}
\end{align}
where $\{e_k\}$ is an arbitrary orthonormal basis for $A^2(\omega)$. The estimate of reproducing kernel of $A^2(\omega)$ has been established by \cite{AH,AP,HLS,LR} for example: for $\omega=e^{-\eta}\in\mathcal{W}$ and $\tau(z)\approx(\Delta\eta(z))^{-1/2}$, there are positive constants $\sigma, C',C''>0$ depending only on $\eta$ such that
\begin{align}\label{kernelesti}
\begin{aligned}
&|K(z,w)|\omega(z)\omega(w)\leq \frac{C'}{\tau(z)\tau(w)}e^{-\sigma d_\tau(z,w)},\qquad\forall z,w\in\D;\\
&|K(z,w)|\omega(z)\omega(w)\geq \frac{C''}{\tau(z)\tau(w)},\qquad d_\tau(z,w)<R;\\
&\|K_z\|_p^p\asymp\omega(z)^{-p}\tau(z)^{-2(p-1)},\quad0<p<\infty.
\end{aligned}
\end{align}

\begin{lem}\label{caltestft}
Given $z, w\in\D$ with $d_\tau(z,w)<R$ for some $R>0$, define
\begin{align*}
f_{z,w}(\xi)=\omega(z)K(\xi,z)(\xi-w)
\end{align*}
then $\|f_{z,w}\|\leq C_2$ where a constant $C_2:=C_2(R,C,C',\sigma)>0$.
\end{lem}
\begin{proof}
Given $0<\delta\leq m_\tau$, we decompose the disk into two parts as follows:
\begin{align*}
\|f_{z,w}\|^2&=\omega(z)^2\int_\D|K(\xi,z)|^2|\xi-w|^2\omega(\xi)^2dA(\xi)\nonumber\\
&=\int_{D(\delta\tau(z))}+\int_{\D\setminus D(\delta\tau(z))}\omega(z)^2|K(\xi,z)|^2|\xi-w|^2\omega(\xi)^2dA(\xi).
\end{align*}
Since $|z-w|<R/C_1\tau(z)$ where $C_1$ is defined in Lemma \ref{setinclus}, (\ref{kernelesti}) and the triangle inequality yield
\begin{align*}
&\omega(z)^2\int_{D(\delta\tau(z))}|K(\xi,z)|^2|\xi-w|^2\omega(\xi)^2dA(\xi)\\
&<\frac{C'}{\tau(z)^2}\int_{D(\delta\tau(z))}\frac{|\xi-z|^2+|z-w|^2}{\tau(\xi)^2}dA(\xi)\lesssim1.
\end{align*}
The constant suppressed depends on $R,C'$ and $\sigma$. On the other hand, we have
\begin{align*}
|z-w|<\frac{R}{C_1}\tau(z)<\frac{R}{\delta C_1}|z-\xi|\quad\text{for}\quad\xi\in\D\setminus D(\delta\tau(z)).
\end{align*}
This, together with (\ref{estimate}) and (\ref{kernelesti}), gives the following inequality:
\begin{align*}
&\omega(z)^2\int_{\D\setminus D(\delta\tau(z))}|K(\xi,z)|^2|\xi-w|^2\omega(\xi)^2dA(\xi)\\
&<2CC'\int_{\D\setminus D(\delta\tau(z))}\frac{|\xi-z|^2+|z-w|^2}{\tau(z)^2\tau(\xi)^2}\left(\frac{\min(\tau(z),\tau(\xi))}{|\xi-z|}\right)^6dA(\xi)\\
&\lesssim\int_{\D\setminus D(\delta\tau(z))}\frac{\tau(z)^2}{|\xi-z|^4}dA(\xi)\leq\sum_{j=0}^{\infty}\int_{2^j\delta\tau(z)<|z-\xi|\leq2^{j+1}\delta\tau(z)}\frac{\tau(z)^2}{|\xi-z|^4}dA(\xi)\\
&\leq\sum_{j=0}^{\infty}\frac{1}{2^{4j}\delta^4\tau(z)^2}2^{2j+2}\delta^2\tau(z)^2\lesssim1.
\end{align*}
All constants suppressed in the above depend on $R,C,C'$ and $\sigma$. It completes the proof.
\end{proof}
Now, we introduce a pseudodistance defined involving the reproducing kernel function:
\begin{align*}
\mathcal{S}(z,w):=\left(1-\frac{|K(z,w)|}{\|K_z\|\|K_w\|}\right)^{1/2},\quad z,w\in\D.
\end{align*}
It is known as Skwarczy\`nski pseudodistance, which is first introduced by M. Skwarczy\`nski in \cite{S}. To compare the distance $\mathcal{S}(z,w)$ and $\rho_\tau(z,w)$, we need the following inequality, which is shown in \cite[Theorem 6.4.3]{JP} with its proof:
\begin{align}\label{lowerds}
\mathcal{S}(z,w)\leq\frac{M(z,w)}{\sqrt{K(z,z)}}\leq\sqrt{2}\mathcal{S}(z,w),\quad z,w\in\D
\end{align}
where $M(z,w):=\sup\{|f(z)|:f\in A^2(\omega), \|f\|=1, f(w)=0\}$. Using this inequality, we obtain the following comparison of $\mathcal{S}(z,w)$ and $\rho_\tau(z,w)$, which is essential to complete our main results.
\begin{thm}\label{key1}
Let $\omega\in\mathcal{W}$. Then there is a constant $C_3:=C_3(C,C',C'',\sigma)>0$ such that $\rho_\tau(z,w)\leq C_3\mathcal{S}(z,w)$ for $z,w\in\D$.
\end{thm}
\begin{proof}
By the kernel estimate (\ref{kernelesti}), there is a constant $C>0$ such that
\begin{align*}
\frac{|K(z,w)|}{\|K_z\|\|K_w\|}\leq Ce^{-\sigma d_\tau(z,w)},\quad z,w\in\D.
\end{align*}
It follows by the inequality above that for $d_\tau(z,w)\geq R$ with $Ce^{-\sigma R}<1/4$,
\begin{align*}
\frac{1}{2}\leq\frac{(1-Ce^{-\sigma d_\tau(z,w)})^{1/2}}{1-e^{-d_\tau(z,w)}}\leq\frac{\mathcal{S}(z,w)}{\rho_\tau(z,w)}.
\end{align*}
Now, let's consider the case of $d_\tau(z,w)<R$ with $\tau(z)\leq\tau(w)$. Applying $f_{z,w}(\xi)=\omega(z)K(\xi,z)(\xi-w)$ in $M(z,w)$ to (\ref{lowerds}), it follows by \eqref{dist} and \eqref{kernelesti} that
\begin{align*}
\mathcal{S}(z,w)\gtrsim\frac{|f_{z,w}(z)|}{\sqrt{K(z,z)}}\geq C''\frac{|z-w|}{\tau(z)}\gtrsim d_\tau(z,w)\geq\rho_\tau(z,w);
\end{align*}
the constants suppressed above depend on $C_2$ appearing in Lemma \ref{caltestft} and $C''$. Therefore, $\rho_\tau(z,w)\lesssim\mathcal{S}(z,w)$ for any $z,w\in\D$. Therefore, we complete our proof.
\end{proof}
Using the result above, we get the following optimal norm estimate for the difference of the reproducing kernels:
\begin{thm}\label{diffkernelnorm}
Let $\omega\in\mathcal{W}$. Then for any $z,w\in\D$,
\begin{align*}
\frac{\|K_z-K_w\|^2}{\|K_z\|^2+\|K_w\|^2}\approx\rho_\tau(z,w)^2.
\end{align*}
\end{thm}
\begin{proof}
By Theorem \ref{key1}, we have
\begin{align*}
\|K_{z}-K_{w}\|^2&\geq\|K_{z}\|^2+\|K_{w}\|^2-2|K(z,w)|\\
&=\|K_{z}\|^2+\|K_{w}\|^2-2\|K_{z}\|\|K_{w}\|(1-\mathcal{S}(z,w)^2)\\
&\geq(\|K_{z}\|^2+\|K_{w}\|^2)\mathcal{S}(z,w)^2\\
&\gtrsim(\|K_{z}\|^2+\|K_{w}\|^2)\rho_\tau(z,w)^2.
\end{align*}
Conversely, we first consider the case $d_\tau(z,w)<R$ where $0<R<\delta^2/12(Ce)^{-\frac{1}{M}}$, $C, M$ appeared in (\ref{estimate}). Then $|z-w|<\delta/6\tau(z)$ by Lemma \ref{setinclus}, so for $|z|\geq|w|$,
\begin{align*}
\|K_{z}-K_{w}\|^2&=\sum_{n=0}^{\infty}\frac{1}{\|z^n\|^2}|z^n-w^n|^2\nonumber\\
&\lesssim\frac{\rho_{\tau}(z,w)^2e^{2\eta(z)}}{\tau(z)^2}
\sum_{n=0}^{\infty}\frac{1}{\|z^n\|^2}\int_{D(\delta\tau(z))}|\xi|^{2n}e^{-2\eta(\xi)}dA(\xi),
\end{align*}
by Lemma \ref{mvtforinterg}. Meanwhile, by (\ref{definition}), (\ref{kernelesti}) and (\ref{equiquan}), we have
\begin{align*}
\sum_{n=0}^{\infty}\frac{1}{\|z^n\|^2}\int_{D(\delta\tau(z))}|\xi|^{2n}e^{-2\eta(\xi)}dA(\xi)
&=\int_{D(\delta\tau(z))}\|K_\xi\|^2e^{-2\eta(\xi)}dA(\xi)\\
&\approx\int_{D(\delta\tau(z))}\frac{1}{\tau(\xi)^2}dA(\xi)\approx1.
\end{align*}
This, together with the kernel estimate (\ref{kernelesti}), yields
\begin{align*}
\|K_{z}-K_{w}\|^2\lesssim\rho_{\tau}(z,w)^2\|K_z\|^2\quad\text{for}\quad|z|\geq|w|.
\end{align*}
In the same way, we have
\begin{align*}
\|K_{z}-K_{w}\|^2\lesssim\rho_{\tau}(z,w)^2\|K_w\|^2\quad\text{for}\quad|z|<|w|.
\end{align*}
This completes the proof for the case $d_\tau(z,w)<R$. For $d_\tau(z,w)\geq R$, we easily obtain from the triangle inequality that
\begin{align*}
\|K_{z}-K_{w}\|^2\lesssim\rho_{\tau}(z,w)^2(\|K_z\|^2+\|K_w\|^2).
\end{align*}
Therefore, we complete our proof.
\end{proof}

\section{Difference of Composition operators}

In this section, we prove our main theorem \ref{k-thm} and theorem \ref{k-thm1}. To prove the compactness of operators, we have the following convenient compactness criterion for a linear combination of weighted composition operators acting on the weighted Bergman spaces.
\begin{lem}\label{cptcriterion}
Let $L$ be a linear combination of composition operators and assume that $L: A^p(\omega)\to A^p(\omega)$ is bounded. Then $L: A^p(\omega)\to A^p(\omega)$ is compact if and only if $\,Lf_k\to 0$ in $A^p(\omega)$ for any bounded sequence $\{f_k\}$ in $A^p(\omega)$ such that $f_k\to 0$ uniformly on compact subsets of $\D$.
\end{lem}
\begin{lem}\cite[Lemma 4.1]{P1}\label{delight}
Let $\omega\in\mathcal{W}$ and $\varphi\in S(\D)$. If there is a curve $\gamma$ connecting to $\zeta\in\partial\D$ and a constant $c>0$ such that $\liminf_{z\to \zeta}\frac{\omega(z)}{\omega(\varphi(z))}\geq c$ where $z\in\gamma$, then
\begin{align*}
\liminf_{z\to \zeta}\frac{\tau(z)}{\tau(\varphi(z))}\geq\min(1,c)\quad\text{for}\quad z\in\gamma.
\end{align*}
Moreover, if $\lim_{z\to \zeta}\frac{\omega(z)}{\omega(\varphi(z))}=0$ then $\limsup_{z\to \zeta}\frac{\tau(z)}{\tau(\varphi(z))}\leq1$.
\end{lem}

\begin{thm}\label{bddnecess}
Let $\omega\in\mathcal{W}$. If $C_\varphi-C_\psi$ is bounded (compact, resp.) on $A^p(\omega)$, then
\begin{align*}
\sup_{z\in\D}\rho_{\tau,\varphi,\psi}(z)^2\left(\frac{\omega(z)}{\omega(\varphi(z))}+\frac{\omega(z)}{\omega(\psi(z))}\right)<\infty \ \  (\rightarrow0\ \ as \ \ |z|\rightarrow1^-, resp.).
\end{align*}
\end{thm}
\begin{proof}
Assume that there exists a sequence $\{z_n\}$ such that
\begin{align}\label{assump1}
\rho_{\tau,\varphi,\psi}(z_n)^2\frac{\omega(z_n)}{\omega(\varphi(z_n))}\rightarrow\infty\quad\text{as}\quad n\rightarrow\infty.
\end{align}
Then we may assume that $\lim_{n\rightarrow\infty}\frac{\omega(z_n)}{\omega(\varphi(z_n))}=\infty$ with $|\varphi(z_n)|\geq|\psi(z_n)|$. Now, define the bounded sequence $\{g_n\}$ on $A^p(\omega)$:
\begin{align}\label{testft}
g_n(\xi)=\frac{K_{\varphi(z_n)}(\xi)}{\omega(\varphi(z_n))^{-1}\tau(\varphi(z_n))^{-2+\frac{2}{p}}}.
\end{align}
By (1) of Lemma \ref{derivsubmean}, Lemma \ref{delight}, Theorem \ref{key1} and \eqref{kernelesti},
\begin{align}\label{goal3}
\|(C_\varphi-C_{\psi})g_{n}\|_p^p&\geq\int_{D(\delta\tau(z_n))}|g_{n}(\varphi(\xi))-g_{n}(\psi(\xi))|^p\omega(\xi)^pdA(\xi)\nonumber\\
&\gtrsim\tau(z_n)^2|g_n(\varphi(z_n))-g_{n}(\psi(z_n))|^p\omega(z_n)^p\nonumber\\
&=|K_{\varphi(z_n)}(\varphi(z_n))-K_{\varphi(z_n)}(\psi(z_n))|^p\frac{\omega(z_n)^p}{\omega(\varphi(z_n))^{-p}}\frac{\tau(z_n)^2}{\tau(\varphi(z_n))^{2(1-p)}}\nonumber\\
&\gtrsim\frac{\omega(z_n)^p}{\omega(\varphi(z_n))^p}\frac{\tau(z_n)^2}{\tau(\varphi(z_n))^{2}}\left|1-\left|
\frac{K_{\varphi(z_n)}(\psi(z_n))}{K_{\varphi(z_n)}(\varphi(z_n))}\right|\right|^p\\
&\gtrsim\frac{\omega(z_n)^p}{\omega(\varphi(z_n))^p}\left|1-\frac{|K_{\varphi(z_n)}(\psi(z_n))|}{\|K_{\varphi(z_n)}\|\|K_{\psi(z_n)}\|}\right|^p\nonumber\\
&\gtrsim\frac{\omega(z_n)^p}{\omega(\varphi(z_n))^p}\rho_{\tau,\varphi,\psi}(z_n)^{2p}\gtrsim1\nonumber
\end{align}
as $n\rightarrow\infty$. This yields a contradiction to (\ref{assump1}) when $n\rightarrow\infty$. The compactness part is immediate from Lemma \ref{cptcriterion} since the sequence $\{g_n\}$ in \eqref{testft} uniformly converges to $0$ on compact subsets of $\D$. This time, we assume that there is a boundary point $\zeta$ and a sequence $\{z_n\}$ such that
\begin{align*}
\lim_{z_n\rightarrow\zeta}\rho_{\tau,\varphi,\psi}(z_n)^2\frac{\omega(z_n)}{\omega(\varphi(z_n))}\neq0.
\end{align*}
Then, $\lim_{n\rightarrow\infty}\frac{\omega(z_n)}{\omega(\varphi(z_n))}\gtrsim1$ with $|\varphi(z_n)|\geq|\psi(z_n)|$ can be assumed. Thus, by the same argument with the boundedness part and Lemma \ref{cptcriterion}, we derive a contradiction to our assumption. This completes all our proof.
\end{proof}
From \eqref{goal3}, we obtain lower bounds of operator norm and the essential operator norm of $C_\varphi-C_\psi$ on $A^p(\omega)$ promptly. In what follows, we denote the operator norm $\|T\|_{A^p(\omega)}$ by $\|T\|_p$ for simplicity.
\begin{cor}\label{necessarybdd}
Let $\omega\in\mathcal{W}$. For $0<p<\infty$,
\begin{align*}
\|C_\varphi-C_{\psi}\|_{p}&\gtrsim\limsup_{\rho_{\tau,\varphi,\psi}(z)\rightarrow1}\left(\frac{\omega(z)}{\omega(\varphi(z))}+
\frac{\omega(z)}{\omega(\psi(z))}\right).
\end{align*}
\end{cor}
\begin{proof}
Note that $|z_n|\rightarrow1$ when $\rho_{\tau,\varphi,\psi}(z_n)\rightarrow1$ as $n\rightarrow\infty$. Then  $\lim_{n\rightarrow\infty}\frac{\tau(z_n)}{\tau(\varphi(z_n))}\gtrsim1$ if $\lim_{n\rightarrow\infty}\frac{\omega(z_n)}{\omega(\varphi(z_n))}\neq0$ by Lemma \ref{delight}. Moreover, $\frac{\tau(z_n)}{\tau(\varphi(z_n))}\lesssim1$ in (\ref{goal3}) when $\lim_{n\rightarrow\infty}\frac{\omega(z_n)}{\omega(\varphi(z_n))}=0$ by Lemma \ref{delight}, thus we obtain our desired lower bound from the proof of Theorem \ref{bddnecess}.
\end{proof}
\begin{cor}
Let $\omega\in\mathcal{W}$. For $0<p<\infty$,
\begin{align*}
\|C_\varphi-C_{\psi}\|_{e}&\gtrsim\limsup_{\rho_{\tau,\varphi,\psi}(z)\rightarrow1}\left(\frac{\omega(z)}{\omega(\varphi(z))}+
\frac{\omega(z)}{\omega(\psi(z))}\right).
\end{align*}
\end{cor}
\begin{proof}
For any compact operator $K$ on $A^p(\omega)$, it follows by Corollary \ref{necessarybdd} that
\begin{align*}
\|C_\varphi-C_{\psi}-K\|_p&\gtrsim\limsup_{n\rightarrow\infty}|\|(C_\varphi-C_{\psi})g_n\|_p-\|Kg_n\|_p|
\\&=\limsup_{n\rightarrow\infty}\|(C_\varphi-C_{\psi})g_n\|_p,
\end{align*}
where $g_n$ is as defined in (\ref{testft}) with any sequence $\{z_n\}$ converging to $1$. The rest of the proof is the same as the proof of Corollary \ref{necessarybdd}.
\end{proof}
As mentioned earlier, not every composition operator is bounded on $A^p(\omega)$, so it is natural to ask of the question when the difference between two composition operators is bounded on $A^p(\omega)$.
\begin{pro}\label{mainresult2}
Let $\omega=e^{-\eta}\in\mathcal{W}$. Then, for $0<p<\infty$,
\begin{align*}
\|C_\varphi-C_\psi\|_{p}\lesssim\|\rho_{\tau,\varphi,\psi} C_\varphi\|_{p}+\|\rho_{\tau,\varphi,\psi} C_\psi\|_{p}.
\end{align*}
\end{pro}
\begin{proof}
For $\|f\|_p\leq1$, we write
\begin{align*}
&\|(C_\varphi-C_\psi)f\|^p_p\\
&=\int_{E}|(f\circ\varphi-f\circ\psi) e^{-\eta}|^pdA+\int_{\D\setminus E}|(f\circ\varphi-f\circ\psi)e^{-\eta}|^pdA
\end{align*}
where
\begin{align}\label{set}
E:=\{z\in\D:d_\tau(\varphi(z),\psi(z))<R\}.
\end{align}
Take $0<R<\delta^2/12(Ce)^{-\frac{1}{M}}$ where $M,C$ in (\ref{estimate}). Then by Lemma \ref{setinclus},
\begin{align*}
|\varphi(z)-\psi(z)|<\delta/6\tau(\varphi(z))\quad\text{for}\quad z\in E.
\end{align*}
Also, we denote $E_1:=\{z\in E:|\varphi(z)|\geq|\psi(z)|\}$. Then by Lemma \ref{mvtforinterg}, we have
\begin{align*}
&\int_{E}|(f\circ\varphi-f\circ\psi)e^{-\eta}|^pdA\\
&\lesssim\int_{E_1}\frac{\rho_{\tau,\varphi,\psi}(z)^pe^{p\eta(\varphi(z))-p\eta(z)}}{\tau(\varphi(z))^2}
\int_{D(\delta\tau(\varphi(z)))}|fe^{-\eta}|^pdA\,dA(z)\\
&+\int_{E\setminus E_1}\frac{\rho_{\tau,\varphi,\psi}(z)^pe^{p\eta(\psi(z))-p\eta(z)}}{\tau(\psi(z))^2}
\int_{D(\delta\tau(\psi(z)))}|fe^{-\eta}|^pdA\,dA(z).
\end{align*}

Therefore, by (\ref{equiquan}), (\ref{opnorm}) and Fubini's theorem, we have
\begin{align}
&\int_{E}|(f\circ\phi-f\circ\psi)e^{-\eta}|^pdA\nonumber\\
&\lesssim\int_{\D}|f(\xi)e^{-\eta(\xi)}|^p\int_{\varphi^{-1}(D(\delta\tau(\xi)))}
\frac{\chi_{E_1}(z)\rho_{\tau,\varphi,\psi}(z)^p}{\tau(\varphi(z))^2}\frac{e^{-p\eta(z)}}{e^{-p\eta(\varphi(z))}}dA(z)dA(\xi)\nonumber\\
&+\int_{\D}|f(\xi)e^{-\eta(\xi)}|^p\int_{\psi^{-1}(D(\delta\tau(\xi)))}
\frac{\chi_{E\setminus E_1}(z)\rho_{\tau,\varphi,\psi}(z)^p}{\tau(\psi(z))^2}\frac{e^{-p\eta(z)}}{e^{-p\eta(\psi(z))}}dA(z)dA(\xi)\nonumber\\
&\lesssim\sup_{\xi\in\D}\frac{\mu_{\rho_{\tau,\varphi,\psi},\varphi,p}(D(\delta\tau(\xi)))}{\tau(\xi)^2}+\sup_{\xi\in\D}\frac{\mu_{\rho_{\tau,\varphi,\psi},\psi,p}(D(\delta\tau(\xi)))}
{\tau(\xi)^2}\nonumber\\
&\lesssim\|\rho_{\tau,\varphi,\psi} C_\varphi\|^p_{p}+\|\rho_{\tau,\varphi,\psi} C_\psi\|^p_{p}.\label{integralE}
\end{align}
Since $\rho_{\tau,\varphi,\psi}(z)\geq1-e^{-R}$ on $\D\setminus E$, we have
\begin{align}\label{integralEc}
\int_{E^c}|f\circ\varphi-f\circ\psi|^pe^{-p\eta} dA&\lesssim\int_{E^c}(|f\circ\varphi|^p+|f\circ\psi|^p)\rho_{\tau,\varphi,\psi}^pe^{-p\eta}dA\nonumber\\
&\leq\|\rho_{\tau,\varphi,\psi} C_\varphi f\|^p_p+\|\rho_{\tau,\varphi,\psi} C_\psi f\|^p_p
\end{align}
for all $f\in A^p(\omega)$. By (\ref{integralE}) and (\ref{integralEc}), we obtain the asserted inequality.
\end{proof}

In the same way as the proof of Proposition \ref{mainresult2}, we can obtain the following sufficient condition for the compactness of $C_\varphi-C_\psi$ immediately.
\begin{pro}\label{mainresult3}
Let $\omega=e^{-\eta}\in\mathcal{W}$ and $0<p<\infty$. If $\rho_{\tau,\varphi,\psi} C_\varphi$ and $\rho_{\tau,\varphi,\psi} C_\psi$ are compact from $A^p(\omega)$ into $L^p(\omega dA)$ then $C_\varphi-C_\psi$ is compact on $A^p(\omega)$.
\end{pro}
\begin{proof}
Consider a sequence $\{f_k\}$ converging to $0$ weakly on $A^p(\omega)$ when $k\rightarrow\infty$. We claim that the following integral vanishes as $k\rightarrow\infty$,
\begin{align}\label{goal}
&\|f_k\circ\varphi-f_k\circ\psi\|_p^p\nonumber\\
&=\int_{E}|f_k\circ\varphi-f_k\circ\psi|^pe^{-p\eta}dA+\int_{\D\setminus E}|f_k\circ\varphi-f_k\circ\psi|^pe^{-p\eta}dA
\end{align}
where the set $E$ is as defined in (\ref{set}). First, we easily see that the second integral of (\ref{goal}) vanishes as $k\rightarrow\infty$ by (\ref{integralEc}) and Lemma \ref{cptcriterion}. Thus, we only need to verify that the first integral of (\ref{goal}) converges to $0$. Since $\rho_{\tau,\varphi,\psi} C_\varphi$ and $\rho_{\tau,\varphi,\psi} C_\psi$ are compact on $A^p(\omega)$, there is $r(\epsilon)>0$ such that for $r<|\xi|<1$,
\begin{align*}
\frac{1}{\tau(\xi)^2}\left(\int_{\varphi^{-1}(D(\delta\tau(\xi)))}
\rho_{\tau,\varphi,\psi}^p\frac{e^{-p\eta}}{e^{-p(\eta\circ\varphi)}}dA+\int_{\psi^{-1}(D(\delta\tau(\xi)))}
\rho_{\tau,\varphi,\psi}^p\frac{e^{-p\eta}}{e^{-p(\eta\circ\psi)}}dA\right)<\epsilon.
\end{align*}
Using the same method in Proposition \ref{mainresult2}, the first integral of (\ref{goal}) is dominated by
\begin{align*}
&\int_{E}|(f_k\circ\varphi-f_k\circ\psi)e^{-\eta}|^pdA\nonumber\\
&\lesssim\int_{\D}|f_k(\xi)|^pe^{-p\eta(\xi)}\int_{\varphi^{-1}(D(\delta\tau(\xi)))}
\frac{\chi_{E_1}(z)\rho_{\tau,\varphi,\psi}(z)^p}{\tau(\varphi(z))^2}\frac{e^{-p\eta(z)}}{e^{-p\eta(\varphi(z))}}dA(z)dA(\xi)\nonumber\\
&+\int_{\D}|f_k(\xi)e^{-\eta(\xi)}|^p\int_{\psi^{-1}(D(\delta\tau(\xi)))}
\frac{\chi_{E\setminus E_1}(z)\rho_{\tau,\varphi,\psi}(z)^p}{\tau(\psi(z))^2}\frac{e^{-p\eta(z)}}{e^{-p\eta(\psi(z))}}dA(z)dA(\xi)\nonumber\\
&\lesssim\epsilon\int_{\D\setminus r\D}|f_k|^pe^{-p\eta}dA+(\|\rho_{\tau,\varphi,\psi} C_\varphi\|_p^p+\|\rho_{\tau,\varphi,\psi} C_\psi\|_p^p)\int_{r\D}|f_k|^pe^{-p\eta}dA.
\end{align*}
Therefore, we can make the integration above small when $k\rightarrow\infty$, so we complete our proof.
\end{proof}
We have shown the implications $(3)\Longrightarrow(1)$ and $(1)\Longrightarrow(2)$ of Theorem \ref{k-thm1} in Proposition \ref{mainresult3} and Theorem \ref{bddnecess}, respectively. Now, we remain to prove the implication $(2)\Longrightarrow(3)$ of Theorem \ref{k-thm1}. The following lemma is an improved version of \cite[Lemma 3.1]{P1}.
\begin{lem}\label{weightedcomp}
Let $\omega\in\mathcal{W}$ and $u$ be a measurable function on $\D$. If $C_\varphi$ be bounded on $A^p(\omega)$ for some $0<p<\infty$ and $s>1$,
\begin{align*}
\lim_{|z|\rightarrow1^-}u(z)^s\frac{\omega(z)}{\omega(\varphi(z))}=0
\end{align*}
then $uC_\varphi:A^p(\omega)\rightarrow L^p(\omega dA)$ is compact for $0<p<\infty$.
\end{lem}
\begin{proof}
For any $\xi\in\varphi^{-1}(D(\delta\tau(z)))$,
\begin{align*}
\delta c_1(1-|z|)>\delta\tau(z)\geq|\varphi(\xi)-z|\geq1-|\varphi(\xi)|-(1-|z|).
\end{align*}
This, together with the Schwartz pick theorem, yields
\begin{align*}
1-|\xi|\lesssim1-|\varphi(\xi)|\lesssim1-|z|.
\end{align*}
For given $\epsilon>0$ and $s>1$, we take $1-|z|<r$ with a sufficiently small $r>0$ so that
\begin{align*}
u(\xi)^s\frac{\omega(\xi)}{\omega(\varphi(\xi))}<\epsilon \ \quad \  on \quad \  \varphi^{-1}(D(\delta\tau(z))).
\end{align*}
This yields
\begin{align*}
\omega^{-p}[u^p\omega^p dA]\circ\varphi^{-1}(D(\delta\tau(z)))&=\int_{\varphi^{-1}(D(\delta\tau(z)))}u(\xi)^p\frac{\omega(\xi)^p}{\omega(\varphi(\xi))^p}\,dA(\xi)\\
&=\int_{\varphi^{-1}(D(\delta\tau(z)))}\left(u(\xi)^{s}\frac{\omega(\xi)}{\omega(\varphi(\xi))}\right)^{\frac{p}{s}}
\left(\frac{\omega(\xi)}{\omega(\varphi(\xi))}\right)^{p\left(1-\frac{1}{s}\right)}\,dA(\xi)\\
&\lesssim\epsilon^{\frac{p}{s}}\int_{\varphi^{-1}(D(\delta\tau(z)))}\left(\frac{\omega(\xi)}{\omega(\varphi(\xi))}\right)^{p\left(1-\frac{1}{s}\right)}\,dA(\xi)\\
&\lesssim\epsilon^{\frac{p}{s}}\tau(z)^2.
\end{align*}
The last inequality is from Theorem \ref{embedding} since $C_\varphi$ is bounded on $A^p(\omega)$ for all $0<p<\infty$ eventually by Remark \ref{bddp}. This, together with the Carleson measure theorem, completes our proof.
\end{proof}
\begin{rmk}
If the measurable function $u$ is bounded on $\D$ in Lemma \ref{weightedcomp}, we can see that $s=1$ can be included in the proof. For details, you can refer to the proof in \cite[Lemma 3.1]{P1}. So, applying $u=\rho_{\tau,\varphi,\psi}$, we easily obtain the following result.
\end{rmk}
\begin{cor}\label{result1}
Let $\omega\in\mathcal{W}$ and $0<p<\infty$. Suppose $C_\varphi, C_\psi$ are bounded on $A^p(\omega)$. If for $s\geq1$
\begin{align*}
\lim_{|z|\rightarrow1}\rho_{\tau,\varphi,\psi}(z)^s\left(\frac{\omega(z)}{\omega(\varphi(z))}+\frac{\omega(z)}{\omega(\psi(z))}\right)=0
\end{align*}
then $\rho_{\tau,\varphi,\psi} C_\varphi$ and $\rho_{\tau,\varphi,\psi} C_\psi$ are compact from $A^p(\omega)$ into $L^p(\omega dA)$.
\end{cor}
Combining Proposition \ref{mainresult3}, Theorem \ref{bddnecess} and Corollary \ref{result1}, we complete the proof of Theorem \ref{k-thm1}.

\section{Example}
Recall that $\varphi$ has a finite angular derivative $\varphi'(\zeta)$ at a boundary point $\zeta$ if we denote $\varphi(\zeta):=\angle\lim_{z\rightarrow\zeta}\varphi(z)$,
\begin{align*}
\varphi'(\zeta):=\angle\lim_{\substack{z\to \zeta\\ z\in \Gamma(\zeta,\alpha)}}\frac{\varphi(z)-\varphi(\zeta)}{z-\zeta}<\infty\quad\text{for each}\quad\alpha>1
\end{align*}
where $\Gamma(\zeta,\alpha)=\{z\in\D:|z-\zeta|<\alpha(1-|z|)\}$. It is well-known as the Julia-Caratheodory Theorem that
\begin{align}\label{Julia}
|\varphi'(\zeta)|=\liminf_{z\rightarrow\zeta}\frac{1-|\varphi(z)|}{1-|z|}<\infty.
\end{align}

In \cite{KM}, Kriete and MacCluer gave the characterization for the boundedness and the compactness of $C_\varphi$ with respect to the angular derivative of $\varphi$ as follows: $C_\varphi$ is unbounded on $A^2(\omega)$ if there exists a boundary point $\zeta$ such that $|\varphi'(\zeta)|<1$. Moreover,
\begin{align}\label{angulcomp}
C_\varphi\quad \text{is compact on } A^2(\omega) \quad \iff \quad |\varphi'(\zeta)|>1,\quad \forall\zeta\in\partial\D.
\end{align}
The following Lemma shows the relation between the value of $\frac{\omega(z)}{\omega(\varphi(z))}$ near the boundary point $\zeta$ and $|\varphi'(\zeta)|$.
\begin{lem}\cite[Lemma 2.9]{P1}\label{pointesti}
Let $\omega\in\mathcal{W}$. If $\lim_{z\rightarrow\zeta}\frac{\omega(z)}{\omega(\varphi(z))}=0$ for $\zeta\in\partial\D$ then $|\varphi'(\zeta)|>1$.
\end{lem}
Consider two analytic self-maps $\varphi$ and $\psi$ having finite
angular derivatives at $\zeta$. We say that $\varphi$ and $\psi$ have the same $M$ order data at $\zeta$ if $\varphi, \psi$ are $M$-th continuously differentiable at $\zeta$ and $\varphi^{(n)}(\zeta)=\psi^{(n)}(\zeta)$ for $n=0,1,\ldots,M$.\\
\indent For $k>0$, it is said that $\varphi(\D)$ has order of contact at most $k$ with $\partial\D$ if for each $\zeta\in\partial\D$ there is a neighborhood $\mathcal{N}(\zeta)\cap\D$ centered at $\zeta$ so that
\begin{align*}
\inf\left\{\frac{1-|\varphi(z)|}{|\varphi(\zeta)-\varphi(z)|^k}:z\in\mathcal{N}(\zeta)\cap\D\right\}>0.
\end{align*}

\begin{pro}\label{omegatau}
Let $\omega=e^{-\eta}\in\mathcal{W}$ and $\varphi, \psi\in\mathcal{S}(\D)$. Suppose there exists the smallest integer $m\geq1$ satisfying $\frac{\tau(r)}{(1-r)^m}\gtrsim1$ and $\varphi, \psi$ are $M$-th continuously differentiable at $\zeta$ for $M\geq1$. If $\varphi, \psi$ have the same $M$ order data at $\zeta$ and $\varphi$ has order of contact at most $\frac{M}{m}$ at $\zeta$ then $\lim_{z\rightarrow\zeta}\rho_{\tau,\varphi,\psi}(z)=0$.
\end{pro}
\begin{proof}
Assume that $\varphi(\zeta)=\psi(\zeta)$, $\varphi^{(n)}(\zeta)=\psi^{(n)}(\zeta)$, $n=1,\ldots,M$, then by the Taylor expansions of $\varphi$ and $\psi$, we have
\begin{align*}
\varphi(z)-\psi(z)=h(z)
\end{align*}
where $h(z)=o(|z-\zeta|^{M})$. Since $\tau(\varphi(z))\gtrsim(1-|\varphi(z))^m$ and $\varphi$ has order of contact $\frac{M}{m}$ at $\zeta$, we obtain
\begin{align*}
\frac{|\varphi(z)-\psi(z)|}{\tau(\varphi(z))}&\lesssim\frac{|h(z)|}{|\varphi(z)-\varphi(\zeta)|^{M}}
\frac{|\varphi(z)-\varphi(\zeta)|^{M}}{(1-|\varphi(z)|)^{m}}\\
&\lesssim\frac{|h(z)|}{|z-\zeta|^{M}}\longrightarrow0
\end{align*}
when $z\rightarrow\zeta$. Thus, we complete our proof by \eqref{dist}.
\end{proof}
We remark that the case $m=1$ of the result above implies the standard weight case, which was shown in \cite{M}. In conjunction with Theorem \ref{k-thm1}, the following result can be obtained immediately from Proposition \ref{omegatau}.

\begin{cor}\label{deriversion}
Let $\omega=e^{-\eta}\in\mathcal{W}$ and $C_\varphi, C_\psi$ be bounded on $A^p(\omega)$. Let
\begin{align*}
F:=\{\zeta\in\partial\D:|\varphi'(\zeta)|=|\psi'(\zeta)|=1\}.
\end{align*}
Suppose there exists the smallest integer $m\geq1$ satisfying $\frac{\tau(r)}{(1-r)^m}\gtrsim1$ and $\varphi, \psi$ are $M$-th continuously differentiable at $\zeta\in F$ for $M\geq1$. If $\varphi, \psi$ have the same $M$ order data at $\zeta\in F$ and $\varphi$ has order of contact at most $\frac{M}{m}$ at $\zeta\in F$ then $C_\varphi-C_\psi$ is compact on $A^p(\omega)$.
\end{cor}
\begin{proof}
By Lemma \ref{pointesti}, if $|\varphi'(\zeta)|\leq1$ then $\lim_{z\rightarrow\zeta}\frac{\omega(z)}{\omega(\varphi(z))}\neq0$. Since $C_\varphi, C_\psi$ are bounded, there is no boundary point $\zeta$ such that $|\varphi'(\zeta)|<1$. Thus, it suffices to show that $\limsup_{z\rightarrow\zeta}\rho_{\tau,\varphi,\psi}(z)=0$ for $\zeta\in F$ by Theorem \ref{k-thm1}. Therefore, it complete the proof by Proposition \ref{omegatau}.
\end{proof}
Finally, we provide an example showing that $C_\varphi$ and $C_\psi$ are not compact on $A^2(\omega)$ but their difference is compact.
\begin{ex}
Consider the Bergman space having the weight $\omega(z)=e^{-\frac{1}{1-|z|}}$. Put
\begin{align*}
\varphi(z)=\frac{1+z^2}{2}\quad\text{and}\quad\psi(z)=\varphi(z)+\epsilon(1-z^2)^5,\quad0<\epsilon<\frac{1}{2^8}.
\end{align*}
Then $C_\varphi$ and $C_\psi$ are not compact on $A^2(\omega)$ but $C_\varphi-C_\psi$ is compact.
\end{ex}
\begin{proof}
First, we note that the following estimates hold:
\begin{align}\label{calculation}
\psi(z)=\varphi(z)+32\epsilon(1-\varphi(z))^5\quad\text{and}\quad1-|\varphi(z)|^2\geq|1-\varphi(z)|^2,\quad z\in\D.
\end{align}
Thus, $\psi\in S(\D)$ since for $0<\epsilon<2^{-8}$,
\begin{align}\label{lesson}
|\psi(z)|\leq|\varphi(z)|+32\epsilon(1-|\varphi(z)|^2)^{2}\leq1.
\end{align}
Here, we note that $|\varphi(\zeta)|=|\psi(\zeta)|=1$ only when $\zeta=1,-1$. Thus, we easily obtain that
\begin{align*}
\lim_{z\rightarrow\zeta\neq1,-1}\frac{\omega(z)}{\omega(\varphi(z))}=\lim_{z\rightarrow\zeta\neq1,-1}\frac{\omega(z)}{\omega(\psi(z))}=0.
\end{align*}
Moreover, $|\varphi'(\zeta)|,|\psi'(\zeta)|>1$ when $\zeta\neq1,-1$ by Lemma \ref{pointesti}. On the other hand, since $|\varphi'(1)|=|\varphi'(-1)|=1=\limsup_{z\rightarrow1,-1}\frac{1-|z|}{1-|\varphi(z)|}$ by (\ref{Julia}) and
\begin{align}\label{calculation2}
\frac{|\varphi(z)|-|z|}{(1-|z|)(1-|\varphi(z)|)}\leq\frac{1}{2}\frac{(1-|z|)^2}{(1-|z|)(1-|\varphi(z)|)}\leq\frac{1}{2}\frac{1-|z|}{1-|\varphi(z)|},
\end{align}
we have
\begin{align*}
\limsup_{z\rightarrow1,-1}\frac{\omega(z)}{\omega(\varphi(z))}&=\limsup_{z\rightarrow1,-1}
\exp\left(\frac{|\varphi(z)|-|z|}{(1-|z|)(1-|\varphi(z)|)}\right)\\&\leq
\limsup_{z\rightarrow1,-1}\exp\left(\frac{1}{2}\frac{1-|z|}{1-|\varphi(z)|}\right)<\infty.
\end{align*}
Likewise, by (\ref{lesson}) and (\ref{calculation2}), we have
\begin{align*}
\frac{|\psi(z)|-|z|}{(1-|z|)(1-|\psi(z)|)}&\leq\frac{|\varphi(z)|-|z|+4(1-|z|+|z|-|\varphi(z)|)^2}{(1-|z|)(1-|\psi(z)|)}\\
&\leq\frac{1/2(1-|z|)^2+8(1-|z|)^2+8(|\varphi(z)|-|z|)^2}{(1-|z|)(1-|\psi(z)|)}\\
&\lesssim\frac{(1-|z|)^2}{(1-|z|)(1-|\psi(z)|)}=\frac{1-|z|}{1-|\psi(z)|}.
\end{align*}
Thus, $\limsup_{z\rightarrow1,-1}\frac{\omega(z)}{\omega(\psi(z))}<\infty$ since $|\psi'(1)|=|\psi'(-1)|=1=\limsup_{z\rightarrow1,-1}\frac{1-|z|}{1-|\psi(z)|}$. Therefore, we have proved that $C_\varphi$ and $C_\psi$ are bounded on $A^2(\omega)$ but they are not compact by (\ref{bddwithoutU}) and (\ref{angulcomp}). Furthermore, it is easily checked that $\varphi, \psi$ have the same $4$-order data at $1,-1$ and $\varphi$ has order of contact at most $2$ from (\ref{calculation}). Letting $\tau(z)=(1-|z|)^{3/2}$, we conclude $C_\varphi-C_\psi$ is compact on $A^2(\omega)$ by Corollary \ref{deriversion}.

\end{proof}

\section{Hilbert-Schmidt Difference}
Let $H_1$ and $H_2$ be Hilbert spaces. Recall that a bounded linear operator $T:H_1\rightarrow H_2$ is Hilbert-Schmidt if
\begin{align*}
\|T\|_{HS}=\sum_{n=0}^{\infty}\|Te_n\|^2<\infty
\end{align*}
where $\{e_n\}$ is an arbitrary orthonormal basis for the Hilbert space $H_1$. Thus, using the definition of the reproducing kernel (\ref{definition}), we easily obtain the Hilbert-Schmidt norm of $uC_\varphi$ on $A^2(\omega)$,
\begin{align}\label{defhs}
\|uC_\varphi\|_{HS}&=\sum_{n=1}^{\infty}\int_{\D}|u(z)|^2|e_n(\varphi(z))|^2\omega(z)^2dA\nonumber\\
&=\int_\D|u(z)|^2\|K_{\varphi(z)}\|^2\omega(z)^2dA(z).
\end{align}

\begin{lem}\label{normHS}
Let $\omega\in\mathcal{W}$ and $u$ be a measurable function on $\D$. Then
\begin{align*}
\left\|u(C_{\varphi}-C_\psi)\right\|^2_{HS}=\int_\D |u(z)|^2\|K_{\varphi(z)}-K_{\psi(z)}\|^2\omega(z)^2dA(z).
\end{align*}
\end{lem}
\begin{proof}
Consider an arbitrary orthonormal basis $\{e_n\}$ for $A^2(\omega)$. By the definition of the Hilbert-Schmidt norm and the reproducing kernel of $A^2(\omega)$, it follows that
\begin{align*}
&\sum_{n=0}^{\infty}\|(uC_\varphi-uC_\psi)e_n\|^2\\&=\sum_{n=0}^{\infty}[\|ue_n(\varphi)\|^2+\|ue_n(\psi)\|^2-2\re\langle ue_n(\varphi),ue_n(\psi)\rangle_\omega]\\
&=\int_\D|u(z)|^2\sum_{n=0}^{\infty}\big[|e_n(\varphi(z))|^2+|e_n(\psi(z))|^2-2\re e_n(\varphi(z))\overline{e_n(\psi(z))}\big]\omega(z)^2dA(z)\\
&=\int_\D|u(z)|^2(\|K_{\varphi(z)}\|^2+\|K_{\psi(z)}\|^2-2\re K(\varphi(z),\psi(z)))\omega(z)^2dA(z)\\
&=\int_\D|u(z)|^2\|K_{\varphi(z)}-K_{\psi(z)}\|^2\omega(z)^2dA(z).
\end{align*}
\end{proof}

Using Lemma \ref{normHS} and Theorem \ref{diffkernelnorm}, we obtain a Hilbert-Schmidt norm estimate for $C_\varphi-C_\psi$ on $A^2(\omega)$ involving the distance $\rho_{\tau,\varphi,\psi}$ as follows.
\begin{pro}\label{key3}
Let $\omega=e^{-\eta}\in\mathcal{W}$. Then
\begin{align*}
\left\|C_{\varphi}-C_\psi\right\|^2_{HS}\approx\int_\D\rho_{\tau,\varphi,\psi}(z)^2(\|K_{\varphi(z)}\|^2+\|K_{\psi(z)}\|^2)
\omega(z)^2dA(z).
\end{align*}
\end{pro}

Now, from (\ref{defhs}) and Proposition \ref{key3}, we obtain the following result promptly.
\begin{thm}\label{mainresult5}
Let $\omega\in\mathcal{W}$. Then $\rho_{\tau,\varphi,\psi} C_\varphi$, $\rho_{\tau,\varphi,\psi} C_\psi$ are Hilbert-Schmidt operators from $A^2(\omega)$ into $L^2(\omega dA)$ if and only if $C_\varphi-C_\psi$ is a Hilbert-Schmidt operator on $A^2(\omega)$.
\end{thm}

Finally, we close this section with observing the path components of $\mathcal{C}_{HS}(A^2(\omega))$, which denotes the space of all composition operators on $A^2(\omega)$ endowed with topology induced by the metric:
\begin{align*}
d(C_\varphi,C_\psi):=\left\{
 \begin{array}{cc}
\frac{\|C_\varphi-C_\psi\|_{HS}}{1+\|C_\varphi-C_\psi\|_{HS}}, &  \|C_\varphi-C_\psi\|_{HS}<\infty\\
1,&  \|C_\varphi-C_\psi\|_{HS}=\infty.\\
 \end{array}
 \right.
\end{align*}
It is said that $C_\varphi$, $C_\psi$ are in the same path component of $\mathcal{C}_{HS}(A^2(\omega))$ if there exists a continuous path $C_{\gamma(s)}$ with respect to $s$ in $\mathcal{C}_{HS}(A^2(\omega))$ such that $\gamma(0)=\varphi$, $\gamma(1)=\psi$. Define the set
\begin{align*}
U(C_\varphi):=\{C_\psi:\|C_\varphi-C_\psi\|_{HS}<\infty\}.
\end{align*}
In the same sense, $U(C_\varphi)$ is said to be a linearly connected component if $C_{(1-s)\varphi+s\psi}$ is continuous with respect to $s$ in $\mathcal{C}_{HS}(A^2(\omega))$ for every $C_\psi\in U(C_\varphi)$. Now, we will prove that the set $U(C_\varphi)$ is a linearly connected component of $C_\varphi$ in $\mathcal{C}_{HS}(A^2(\omega))$. \cite{HJM} and \cite{CHK} contain some results for the connected component of the space of composition operators acting on the Hardy and the standard weighted Bergman spaces under the Hilbert-Schmidt norm topologies, respectively.
\begin{lem}\label{careful}
Let $\omega\in\mathcal{W}$. For $z,w\in\D$, denote $z_s:=(1-s)z+sw$ for $s\in[0,1]$. Then
\begin{align*}
\rho_\tau(z_s,z_t)\leq C\rho_\tau(z,w)
\end{align*}
where $C>0$ is independent of $s,t$.
\end{lem}
\begin{proof}
We first consider the case $d_\tau(z,w)<R$. By Lemma \ref{setinclus}, there exists $0<C_1<1$ satisfying
\begin{align*}
d_\tau(z_s,z_t)&\leq\frac{|z_t-z_s|}{\min(\tau(z_s),\tau(z_t))}\leq\frac{|s-t||z-w|}{\min(\tau(z),\tau(w))}\leq\frac{1}{C_1}d_\tau(z,w).
\end{align*}
Moreover, it is clear that $\rho_\tau(z_s,z_t)\leq \frac{1}{1-e^{-R}}\rho_\tau(z,w)$ for $d_\tau(z,w)\geq R$. This completes the proof.
\end{proof}
\begin{thm}\label{mainresult1}
Let $\omega\in\mathcal{W}$. Then $C_\varphi-C_\psi$ is Hilbert-Schmidt on $A^2(\omega)$ if and only if $C_\varphi$ and $C_\psi$ lie in the same path component of $\mathcal{C}_{HS}(A^2(\omega))$.
\end{thm}
\begin{proof}
Assume that there exists a continuous path $C_{\gamma(s)}:[0,1]\rightarrow \mathcal{C}_{HS}(A^2(\omega))$. Then $C_{\gamma(s)}$ is uniformly continuous on $[0,1]$ so that given $\epsilon>0$, there is a partition $\{s_0=0,s_1,\ldots,s_{N-1},s_N=1\}\subset[0,1]$ such that $\|C_{\gamma(s_i)}-C_{\gamma(s_{i+1})}\|_{HS}<\epsilon$ for all $i=0,1,\ldots,N-1$. Thus $C_\varphi-C_\psi$ is Hilbert-Schmidt on $A^2(\omega)$ by the triangle inequality. To show the necessity, we will show that
\begin{align*}
\lim_{t\rightarrow s}\|C_{\varphi_s}-C_{\varphi_t}\|_{HS}=0
\end{align*}
when $\varphi_s=(1-s)\varphi+s\psi$ for $0\leq s\leq1$. Since $|\varphi_s(z)|\leq\max(|\varphi(z)|,|\psi(z)|)$ and $\|K_z\|$ increases with $|z|$ by (\ref{kernelesti}), we have $\|K_{\varphi_s(z)}\|\lesssim\|K_{\varphi(z)}\|+\|K_{\psi(z)}\|$. Therefore, using the results of Lemma \ref{normHS}, Proposition \ref{key3} and Lemma \ref{careful}, we have
\begin{align}\label{goal4}
\|C_{\varphi_s}-C_{\varphi_t}\|_{HS}&\lesssim\int_{\D}\rho_{\tau,\varphi_s,\varphi_t}(z)^2(\|K_{\varphi_s(z)}\|^2+\|K_{\varphi_t(z)}\|^2)\omega(z)^2dA(z)\\
&\lesssim\int_{\D}\rho_{\tau,\varphi,\psi}(z)^2(\|K_{\varphi(z)}\|^2+\|K_{\psi(z)}\|^2)\omega(z)^2dA(z)<\infty\nonumber.
\end{align}
By the Dominated Convergence Theorem, (\ref{goal4}) vanishes when $t\rightarrow s$ for $\rho_{\tau,\varphi_s,\varphi_t}(z)\rightarrow0$.
\end{proof}
The inequality above gives the following result immediately.
\begin{cor}
Let $\omega\in\mathcal{W}$. Define $\varphi_s=(1-s)\varphi+s\psi$ for $0\leq s\leq1$. Then $C_\varphi-C_\psi$ is Hilbert-Schmidt on $A^2(\omega)$ if and only if $C_{\varphi_s}-C_{\varphi_t}$ is Hilbert-Schmidt on $A^2(\omega)$ for $0\leq s,t\leq1$.
\end{cor}

\bibliographystyle{amsplain}

\end{document}